\documentclass[preprint,fleqn,sort&compress,3p]{elsarticle}
\usepackage{mathrsfs}
\usepackage{amsfonts}

\usepackage{graphicx}
\usepackage{amssymb}
 \usepackage{amsthm,amsmath}
\usepackage{bm}
\usepackage{latexsym}
\usepackage{booktabs}
\usepackage{hyperref}

\newtheorem{Theorem}{Theorem}[section]

\newtheorem{Remark}{Remark}[section]

\usepackage{lipsum}
\makeatletter
\def\ps@pprintTitle{%
 \let\@oddhead\@empty
 \let\@evenhead\@empty
 \def\@oddfoot{}%
 \let\@evenfoot\@oddfoot}
\makeatother

\begin{document}
\title{The Weighted Arithmetic Mean-Geometric Mean Inequality is Equivalent to the H\"{o}lder Inequality}


\author[a]{Yongtao Li}
\ead{yongtao\_li@csu.edu.cn}
\author[b,c]{Xian-Ming Gu\corref{cor1}}
\ead{guxianming@live.cn}
\cortext[cor1]{Corresponding author}
\author[d]{Jianxing Zhao}
\ead{zhangyanz@mst.edu}
\address[a]{School of Mathematics and Statistics, Central South University, Changsha, Hunan 410083,
China}
\address[b]{School of Economic Mathematics/Institute of Mathematics,\\
Southwestern University of Finance and Economics, Chengdu, Sichuan 611130, China}
\address[c]{Bernoulli Institute for Mathematics, Computer Science and Artificial Intelligence,\\ University of
Groningen, Nijenborgh 9, P.O. Box 407, 9700 AK Groningen, The Netherlands}
\address[d]{College of Data Science and Information Engineering,\\
Guizhou Minzu University, Guiyang, Guizhou~550025,~China}

\begin{abstract}
In the current note, we investigate the mathematical relations
among the weighted arithmetic mean--geometric mean (AM--GM) inequality,
the H\"{o}lder inequality and the weighted power-mean inequality.
Meanwhile, the proofs of mathematical equivalence among the weighted
AM--GM inequality, the weighted power-mean inequality and the H\"{o}lder
inequality {are fully achieved. The new results
are more generalized than those of previous studies.}

\vspace{4mm}
\end{abstract}


\maketitle

\section{Introduction}
In the field of classical analysis, the weighted {arithmetic mean--geometric
mean (AM--GM)} inequality (see e.g., \cite{ZCIT}, pp. 74--75) is often inferred from Jensen's
inequality, which is a more generalized inequality {compared to} the AM--GM inequality;
refer to, e.g., \cite{GJG34,ZCIT}. In addition, the~H\"{o}lder inequality~\cite{GJG34}
found by \textit{Leonard James Rogers} (1888) and discovered independently by \textit{Otto H\"{o}lder}
(1889) is a basic and indispensable inequality for studying integrals and $L^p$ spaces, and is also an extension
of the {Cauchy--Bunyakovsky--Schwarz (CBS) inequality} \cite{SSDA}. The H\"{o}lder inequality is used to prove the Minkowski inequality,
which is the triangle inequality in the space $L^p(\mu)$ \cite{EFR61,DS70}. The~weighted power mean
(also known as the generalized mean) $M_r^m(a)$ for a sequence $a=(a_1,a_2,\ldots,a_n)$ is defined as
$M_r^m(a)=(m_1a_1^r+m_2a_2^r+\cdots+m_na_n^r)^{\frac{1}{r}}$, which is a family of functions for aggregating
sets of numbers, and plays a vital role in analytical inequalities; see \cite{SBJ97,GJG34} for instance.

In recent years, many researchers have been interested in studying the mathematical equivalence among some
famous analytical inequalities, such as the Cauchy--Schwarz inequality, the Bernoulli inequality, the
Wielandt inequality, and the Minkowski inequality; see \cite{FFNT97,LS06,Ma01,Ma12,Zh01,AIB11,SMSG} for
details{. Additionally, these studies note the relations among the weighted AM--GM inequality, the
H\"{o}lder inequality, and the weighted power-mean inequality are still less clear, although one inequality is often helpful to prove another inequality \cite{ZCIT,AIB11}}. Motivated by these aforementioned studies, in the present note, the~mathematical equivalence
among three such well-known inequalities is proved in detail; the result introduced in \cite{Lin12} is also extended.

The rest of the present note is organized as follows. In the next section, we will present the detailed proofs
of mathematical equivalence among three celebrated mathematical inequalities. Finally, the paper ends with
several concluding remarks in Section \ref{sec2}.

\label{sec1}
{The} weighted AM--GM inequality, the H\"{o}lder inequality,
and the weighted power-mean inequality~{\cite{ZCIT}} (pp. 111--112, Theorem 10.5) are first reviewed {
and they are often related to each other.} Then the results of mathematical
equivalence among three such inequalities will be shown.

\noindent{\bf Weighted AM--GM Inequality.}  If $0\leq c_i\in\mathbb{R}~(i = 1,\ldots,n)$ and $0\leq \lambda_i\in\mathbb{R}~(i = 1,\ldots,n)$
such that $\sum\limits^{n}_{i=1}\lambda_i=1$, then
\begin{equation}
\prod\limits_{k=1}^n c_k^{\lambda_k} \le \sum\limits_{k=1}^n \lambda_kc_k.
\label{e1}
\end{equation}

\noindent{\bf H\"{o}lder Inequality.}  If $0\leq a_i,b_i\in\mathbb{R}~(i = 1,\ldots,n)$ and $p,q\in\mathbb{R}^{+}$ such that $p^{-1}+q^{-1}$ = 1, then
\begin{equation}
\sum_{k=1}^n a_kb_k \le \left( \sum\limits_{k=1}^n a_k^{p}\right)^{\frac{1}{p}}\left( \sum\limits_{k=1}^n b_k^{q}\right)^{\frac{1}{q}}.
\label{e2}
\end{equation}

\noindent{\bf Weighted Power-Mean Inequality.} If $0\leq c_i,~\lambda_i\in\mathbb{R}~(i = 1,\ldots,n)$
such that $\sum\limits^{n}_{i=1}\lambda_i=1$, and $r,s\in\mathbb{R}^{+}$ such that $r \le s$, then
\begin{equation}
\left(\sum\limits_{k=1}^n \lambda_kc_k^r \right)^{\frac{1}{r}} \le \left(\sum\limits_{k=1}^n \lambda_kc_k^s \right)^{\frac{1}{s}}.
\label{e3}
\end{equation}

The word ``equivalence" between two statements $\bold A$ and $\bold B$, by
convention, is understood as follows: $\bold A$ implies $\bold B$ and $\bold B$ implies
$\bold A$. {Two equivalent sentences have the same truth value.} Thus,~this~note
reveals a connection (in the sense of art) between these two well-known facts.
\begin{Theorem}
The H\"{o}lder inequality is equivalent to the weighted AM--GM inequality.
\label{the1}
\end{Theorem}
\begin{proof}  To show that (\ref{e2}) implies  (\ref{e1}), let $a_k=(\lambda_kc_k)^{\frac{1}{p}},
b_k=(\lambda_k)^{\frac{1}{q}}$ in (\ref{e2}) for all $k$; then
\begin{equation}
\sum\limits_{k=1}^n \left[(\lambda_kc_k)^{\frac{1}{p}}(\lambda_k)^{\frac{1}{q}}\right] \le \left( \sum\limits_{k=1}^n \lambda_kc_k \right)^{\frac{1}{p}}\left( \sum\limits_{k=1}^n \lambda_k \right)^{\frac{1}{q}}.
\label{guadd1}
\end{equation}

Since $p^{-1}+q^{-1}$ = 1 and $\lambda_1+\lambda_2+\cdots+\lambda_n=1${,
the inequality (\ref{guadd1}) can be rewritten as:}
\begin{equation}
\sum\limits_{k=1}^n \lambda_kc_k \ge \left(\sum\limits_{k=1}^n \lambda_kc_k^{\frac{1}{p}} \right)^{p}.
\label{e4}
\end{equation}

Now using the inequality (\ref{e4}) successively, it follows that
\begin{equation}
\begin{split}
\sum\limits_{k=1}^n \lambda_kc_k &\ge \left(\sum\limits_{k=1}^n \lambda_kc_k^{\frac{1}{p}} \right)^{p}\ge \left(\sum\limits_{k=1}^n \lambda_kc_k^{\frac{1}{p^{2}}} \right)^{p^{2}} \\ &\ge \cdots \ge \left(\sum\limits_{k=1}^n \lambda_kc_k^{\frac{1}{p^{m}}} \right)^{p^{m}} \ge \cdots.
\end{split}
\label{e5}
\end{equation}

By L'Hospital's rule, it is easy to see that
\begin{equation*}
\lim\limits_{x\rightarrow 0^+} \ln \left(\sum\limits_{k=1}^n \lambda_kc_k^x \right)\Big/x = \sum\limits_{k=1}^n \lambda_k \ln c_k.
\end{equation*}

Thus,
\begin{equation*}
\lim \limits_{x\rightarrow 0^+} \left(\sum\limits_{k=1}^n \lambda_kc_k^x \right)^{\frac{1}{x}} = \prod \limits_{k=1}^n c_k^{\lambda_k}.
\end{equation*}

Thus, in (\ref{e5}), we can pass to the limit by $m \rightarrow + \infty $, giving $ \sum\limits_{k=1}^n \lambda_kc_k \ge \prod\limits_{k=1}^n c_k^{\lambda_k}$; hence, (\ref{e2}) implies (\ref{e1}).

To show the converse, all that is needed is a special case of (\ref{e1}),
\begin{equation}
\lambda_1c_1+\lambda_2c_2 \ge c_1^{\lambda_1}c_2^{\lambda_2}.
\label{e6}
\end{equation}

Since $p^{-1}+q^{-1}$=1, and by (\ref{e6}), thus
\begin{equation*}
\frac{1}{p}a_k^p\sum\limits_{k=1}^n b_k^q + \frac{1}{q}b_k^q\sum\limits_{k=1}^n a_k^p \ge a_kb_k \left(\sum\limits_{k=1}^n b_k^q \right)^{\frac{1}{p}}\left(\sum\limits_{k=1}^n a_k^p \right)^{\frac{1}{q}}.
\end{equation*}

Summing over $k=1,2,\ldots,n$,~it follows that
\begin{equation*}
\sum\limits_{k=1}^n a_k^p \sum\limits_{k=1}^n b_k^q \ge \sum\limits_{k=1}^n
a_kb_k \left( \sum\limits_{k=1}^n b_k^{q}\right)^{\frac{1}{p}}\left( \sum\limits_{k=1}^n a_k^{p}\right)^{\frac{1}{q}}.
\end{equation*}

Thus, (\ref{e1}) implies (\ref{e2}).
\end{proof}

\begin{Remark}
The CBS inequality (the AM--GM inequality) is the special case of the H\"{o}lder inequality (the~weighted AM--GM inequality);
therefore, Theorem \ref{the1} is a generalization of the result established in \cite{Lin12}.
\end{Remark}

\begin{Theorem}
The H\"{o}lder inequality is equivalent to the weighted power-mean inequality.
\end{Theorem}
\begin{proof}
For $r\le s$, let $p=\frac{s}{r}\ge 1$ and $c_k=d^{s}_k$ in (\ref{e4}) for $k = 1,\ldots,n$, then
\begin{equation}
\sum\limits_{k=1}^n \lambda_kd^{s}_k \ge \left(\sum\limits_{k=1}^n \lambda_kd_k^r \right)^{\frac{s}{r}}.
\label{e7x}
\end{equation}

Thus, rewriting (\ref{e7x}) as $\left(\sum\limits_{k=1}^n \lambda_kd_k^s \right)^{\frac{1}{s}}
\ge \left(\sum\limits_{k=1}^n \lambda_kd_k^r \right)^{\frac{1}{r}}$.
Now the task is to prove that (\ref{e3}) implies (\ref{e2});
let $\lambda_k=b_k^{q}\bigl( \sum\limits_{k=1}^n b_k^q\bigr)^{-1},
d_k=a_k/b_k^{q-1},r=1,s=p$ in (\ref{e3}), then
\begin{equation*}
\begin{split}
\sum\limits_{k=1}^na_kb_k
&=\left( \sum\limits_{k=1}^nb_k^q\right)\left( \sum\limits_{k=1}^n
\frac{b_k^q}{\sum\limits_{k=1}^nb_k^q}\cdot \frac{a_k}{b_k^{q-1}}\right) \\
&\le \left( \sum\limits_{k=1}^nb_k^q\right)\left( \sum\limits_{k=1}^n
\frac{b_k^q}{\sum\limits_{k=1}^nb_k^q}\cdot \left(\frac{a_k}{b_k^{q-1}}\right)^p\right)^{\frac{1}{p}} \\
&= \left( \sum\limits_{k=1}^n a_k^{p}\right)^{\frac{1}{p}}
\left( \sum\limits_{k=1}^n b_k^{q}\right)^{\frac{1}{q}}. \quad\qedhere
\end{split}
\end{equation*}
\end{proof}
\begin{Remark}\label{the2}
Here, we can give another proof that (\ref{e3}) implies (\ref{e2}).
Repeatedly using inequality (\ref{e3}), it follows that
\newpage
\begin{equation*}
\sum\limits_{k=1}^n \lambda_kc_k
\ge \left(\sum\limits_{k=1}^n \lambda_kc_k^{\frac{1}{2}} \right)^2
\ge \left(\sum\limits_{k=1}^n \lambda_kc_k^{\frac{1}{3}} \right)^3
\ge \cdots
\ge \left(\sum\limits_{k=1}^n \lambda_kc_k^{\frac{1}{n}} \right)^n
\ge \cdots.
\label{e8}
\end{equation*}

By L'Hospital's rule, it is easy to see that
\begin{equation*}
\sum\limits_{k=1}^n \lambda_kc_k
\ge \prod\limits_{k=1}^n c_k^{\lambda_k}.
\label{e9}
\end{equation*}

By using analogous methods from Theorem \ref{the1}, it can be proved that
\begin{equation*}
\sum\limits_{k=1}^n a_k^p \sum\limits_{k=1}^n b_k^q
\ge \sum\limits_{k=1}^n
a_kb_k \left( \sum\limits_{k=1}^n b_k^{q}\right)^{\frac{1}{p}}
\left( \sum\limits_{k=1}^n a_k^{p}\right)^{\frac{1}{q}}.
\end{equation*}
\end{Remark}
\begin{Theorem}
The weighted power-mean inequality is equivalent to the weighted AM--GM inequality.
\label{the3}
\end{Theorem}
\begin{proof}  To show that (\ref{e1}) implies (\ref{e3}), we merely {exploit} a special case of (\ref{e1}),
\begin{equation}
a_1^{\lambda_1}a_2^{\lambda_2} \le \lambda_1a_1+\lambda_2a_2.
\label{e10}
\end{equation}

Here, we define $U_n(a)=\lambda_1a_1^s+\lambda_2a_2^s+\cdots+\lambda_na_n^s$;
let $a_1=\lambda_ka_k^s (U_n(a))^{-1}, a_2= \lambda_k $ and $\lambda_1=\frac{r}{s},\lambda_2=1-\frac{r}{s}$ in (\ref{e10}), then
\begin{equation*}
\lambda_ka_k^r (U_n(a))^{-\frac{r}{s}}  \le \frac{r}{s} \cdot \lambda_ka_k^s (U_n(a))^{-1}+\left(1-\frac{r}{s}\right)\cdot\lambda_k.
\end{equation*}

Summing over $k=1,2,\ldots,n$, then
\begin{equation*}
\sum\limits_{k=1}^n\lambda_ka_k^r (U_n(a))^{-\frac{r}{s}} \le \sum\limits_{k=1}^n\left[\frac{r}{s} \cdot \lambda_ka_k^s (U_n(a))^{-1}+\left(1-\frac{r}{s}\right)\cdot\lambda_k \right]=1.
\end{equation*}

Therefore,
\begin{equation*}
\left(\sum\limits_{k=1}^n \lambda_kc_k^r \right)^{\frac{1}{r}} \le \left(\sum\limits_{k=1}^n \lambda_kc_k^s \right)^{\frac{1}{s}}.
\end{equation*}

The converse is trivial from Remark \ref{the2}.
\end{proof}
\section{Concluding Remarks}
\label{sec2}
In this note, the mathematical equivalence among the weighted AM--GM inequality,
the H\"{o}lder inequality, and the weighted power-mean inequality is investigated
in detail. Moreover, the interesting conclusions of Lin's paper \cite{Lin12} are
also {extended.} At the end of the present study{, for convenience,
the}~results on the equivalence of some well-known analytical inequalities can
be summarized as~follows:
\begin{itemize}
\item Equivalence of the H\"{o}lder's inequality and the Minkowski inequality; see \cite{Ma01}.
\item Equivalence of the Cauchy--Schwarz inequality and the H\"{o}lder's inequality; see \cite{LS06}.
\item Equivalence of the Cauchy--Schwarz inequality and the Covariance--Variance inequality; see \cite{FFNT97}.
\item Equivalence of the Kantorovich inequality and the Wielandt inequality; see e.g., \cite{Zh01}
\item Equivalence of the AM--GM inequality and the Bernoulli inequality; see e.g., \cite{Ma12}
\item Equivalence of the H\"{o}lder inequality and Artin's theorem; {see} e.g., \cite{AIB11} (pp. 657--663) for
details.
\item Equivalence of the H\"{o}lder inequality and the weighted AM--GM inequality; refer to Theorem \ref{the1}.
\item Equivalence of the H\"{o}lder inequality and the weighted power-mean inequality; {see} Theorem \ref{the2}.
\item Equivalence of the weighted power-mean inequality and the weighted AM--GM inequality; refer to Theorem \ref{the3}.
\end{itemize}

\section*{Author Contributions}
All the authors inferred the main conclusions and approved the current version of this~manuscript.


%

\section*{Acknowledgments}
{\em The authors would like to thank Dr. Minghua Lin for his valuable discussions, which~considerably improved the presentation of our manuscript. This work was funded by the National Natural Science Foundation of China (Grant No. 11501141), Science and Technology Top-notch Talents Support Project of Education Department of Guizhou Province (Grant~No.~QJHKYZ [2016]066).}

\section*{Abbreviations}
\noindent The following abbreviations are used in this manuscript:\\

\noindent
\begin{tabular}{@{}ll}
AM--GM & Arithmetic Mean--Geometric Mean\\
CBS   & Cauchy--Bunyakovsky--Schwarz
\end{tabular}

\vspace{5mm}


\end{document}